\documentclass[a4paper]{article}

\usepackage[english]{babel}
\usepackage[utf8x]{inputenc}
\usepackage{amsmath,amssymb,amstext,amsthm}
\usepackage{graphicx}
\usepackage[colorinlistoftodos]{todonotes}

\usepackage{hyperref}
\usepackage{mathtools}

\newtheorem{theorem}{Theorem}[section]
\newtheorem{lemma}[theorem]{Lemma}
\newtheorem{corollary}[theorem]{Corollary}

\theoremstyle{definition}
\newtheorem*{definition}{Definition}

\theoremstyle{definition}
\newtheorem*{remark}{Remark}

\newcommand{\qarrow}{\xrightarrow{q}}
\newcommand{\earrow}{\xrightarrow{*}}
\newcommand{\tob}{\xrightarrow{B}}
\newcommand{\tov}{\xrightarrow{V}}
\newcommand{\toplus}{\xrightarrow{+}}

\newcommand{\bartheta}{\bar{\vartheta}}
\newcommand{\barthetap}{\bar{\vartheta}^+}
\newcommand{\barthetam}{\bar{\vartheta}^-}
\newcommand{\varthetap}{\vartheta^+}
\newcommand{\varthetam}{\vartheta^-}

\DeclareMathOperator{\tr}{Tr}

\DeclareMathOperator{\vect}{vec}
\DeclareMathOperator{\rk}{rk}

\DeclareMathOperator{\aut}{Aut}

\newcommand{\K}{\mathcal{K}}
\newcommand{\thetak}{\vartheta^\K}
\newcommand{\Thetak}{\Theta^\K}
\newcommand{\barthetak}{\bar{\vartheta}^\K}
\newcommand{\cp}{\mathcal{CP}}
\newcommand{\cpsd}{{\mathcal{CS}_+}}
\newcommand{\dnn}{\mathcal{DNN}}
\newcommand{\psd}{{\mathcal{S}_+}}

\newcommand{\psdn}{{\mathcal{S}^n_+}}
\newcommand{\tok}{\xrightarrow{\K}}
\newcommand{\tokk}{\xrightarrow[\K]{}}
\newcommand{\tocpt}{\xrightarrow{\cp}}
\newcommand{\tocpb}{\xrightarrow[\cp]{}}
\newcommand{\tocpsdt}{\xrightarrow{\cpsd}}
\newcommand{\tocpsdb}{\xrightarrow[\cpsd]{}}
\newcommand{\todnnt}{\xrightarrow{\dnn}}
\newcommand{\todnnb}{\xrightarrow[\dnn]{}}
\newcommand{\topsdt}{\xrightarrow{\psd}}

\newcommand{\cones}{\{\cp, \cpsd, \dnn, \psd\}}
\newcommand{\conesm}{\{\cp, \dnn, \psd\}}
\newcommand{\cone}{\cp, \ \cpsd, \ \dnn, \text{ and } \psd}
\newcommand{\thetacpsd}{\vartheta^\cpsd}
\newcommand{\Thetacpsd}{\Theta^\cpsd}

\newcommand{\qedn}{\qed\vspace{.2cm}}

\title{Conic Formulations of Graph Homomorphisms}
\author{David E. Roberson}

\begin{document}
\maketitle

\begin{abstract}
Given graphs $X$ and $Y$, we define two conic feasibility programs which we show have a solution over the completely positive cone if and only if there exists a homomorphism from $X$ to $Y$. By varying the cone, we obtain similar characterizations of quantum/entanglement-assisted homomorphisms and three previously studied relaxations of these relations. Motivated by this, we investigate the properties of these ``conic homomorphisms" for general (suitable) cones. We also consider two generalized versions of the Lov\'{a}sz theta function, and how they interact with these conic homomorphisms. We prove analogs of several results on classical graph homomorphisms as well as some monotonicity theorems. We also show that one of the generalized theta functions is multiplicative on lexicographic and disjunctive graph products.
\end{abstract}

\section{Introduction}

A map $\varphi$ from the vertex set of a graph $X$ to that of a graph $Y$ is a \emph{homomorphism} if $\varphi(x) \sim \varphi(x')$ whenever $x \sim x'$, where ``$\sim$" denotes adjacency. If such a map exists we say that $X$ has a homomorphism to $Y$ and write $X \to Y$. In this paper we will investigate homomorphisms and relaxations of homomorphisms through the use of conic feasibility programs over various cones. In doing so, we will see that this approach unifies a great many results and ideas that previously seemed distinct from one another.

In~\cite{Piovesan14}, a new approach to quantum/entanglement-assisted chromatic and independence numbers was introduced. Laurent and Piovesan identified the completely positive semidefinite cone and showed that these parameters could be written as optimization programs over this cone. As an intermediate step, they gave a conic feasibility programs which are equivalent to deciding if these parameters at most/least some integer. For the classical version of these parameters, such a decision is equivalent to deciding whether a homomorphism to/from a certain complete graph existed. With this in mind, it becomes natural to ask whether their approach can be extended to quantum/entanglement-assisted homomorphisms in general. We will see that this approach can indeed be extended, and we will give conic feasibility problem formulations for both quantum and entanglement-assisted homomorphisms.

Though the above quantum concerns were the main motivation for this investigation, we found that by varying the cone in question we could obtain various relaxations of quantum homomorphisms that were investigated in~\cite{lovaszhomo,bacik95,Stahlke14}. By considering the completely positive cone we are even able to recover the usual graph homomorphism. The consideration of these other cones was in fact motivated by~\cite{Piovesan14} as well, since there the authors investigate their feasibility programs over the completely positive and doubly nonnegative cones. The main cones we will focus on here are the completely positive ($\cp$), completely positive semidefinite ($\cpsd$), doubly nonnegative ($\dnn$), and positive semidefinite cones ($\psd$). By considering these we will be able to recover very many known results, as well as new ones. Most importantly though, our approach will allow us to prove these results for all of these cones simultaneously, giving us a unified view of things.

After seeing how this conic approach can be used to bring together many results under a single framework, we began to develop a general theory of homomorphisms over cones. Specifically, for any cone which has a few key properties outlined below, we define two types of homomorphisms over this cone. We refer to such homomorphisms as \emph{conic homomorphisms}. For the completely positive semidefinite cone these two definitions correspond to quantum and entanglement-assisted homomorphisms. Interestingly, for each of the completely positive and doubly nonnegative cones the two types of homomorphisms are equivalent. For the doubly nonnegative cone we will see that the two types of homomorphisms we obtain are equivalent to the relations $\tov$ and $\toplus$ from~\cite{Stahlke14}. This implies that $\tov$ and $\toplus$ are actually equivalent, simplifying the hierarchy of relaxations of quantum homomorphisms given in~\cite{Stahlke14}.

In~\cite{chif}, Dukanovic and Rendl introduce two generalized versions of the Lov\'{a}sz theta function. These are defined by considering primal and dual semidefinite programs for Lo\'{a}sz theta and optimizing these over some suitably chosen cone instead of the positive semidefinite cone. These two generalizations will turn out to be intricately connected to our conic homomorphisms, and by considering both notions together we will recover a substantial portion of the known bounds and monotonicity theorems concerning quantum/entanglement-assisted homomorphisms.

\section{Quantum and Entanglement-Assisted Homomorphisms}

Though our main focus will be a purely mathematical analysis of general conic homomorphisms, since our initial motivation came from the study of quantum and entanglement-assisted homomorphisms, we will give a brief introduction to the physical scenarios behind these mathematical notions.

\subsection{Quantum Homomorphisms}

The homomorphism game was originally defined in~\cite{Roberson12} in order to generalize the quantum chromatic number to quantum homomorphisms. For graphs $X$ and $Y$, the $(X,Y)$-homomorphism game is played as follows: a referee sends each of Alice and Bob a vertex of $X$, and they each respond with a vertex of $Y$. In order to win, Alice and Bob must respond with identical vertices if the vertices they received were the same, and they must respond with adjacent vertices if they were sent adjacent vertices. Alice and Bob are allowed to agree on a strategy beforehand, but they are unable to communicate during the game. It is not too difficult to see that, even if given access to some shared randomness, classical players can win if and only if there exists a homomorphism from $X$ to $Y$. However, if players are allowed to make quantum measurements on a shared entangled state, they can sometimes win perfectly even in the absence of such a homomorphism~\cite{BCW98, Galliard02, Avis}. This motivated the definition of quantum homomorphism, which is said to exist from $X$ to $Y$ whenever there is a perfect quantum strategy for the $(X,Y)$-homomorphism game. Mathematically we have the following~\cite{Roberson12}:

\begin{definition}\label{def:qhomo}
We say that there exists a \emph{quantum homomorphism} from a graph $X$ to a graph $Y$, and write $X\qarrow Y$, if there exist positive semidefinite matrices $\rho_{xy}$ for $x \in V(X), y \in V(Y)$ and $\rho$ satisfying
\begin{align*}
&\langle \rho, \rho \rangle = 1 \\
&\sum_{y \in V(Y)} \rho_{xy} = \rho \text{ for all } x \in V(X) \\
&\langle \rho_{xy}, \rho_{x'y'} \rangle = 0 \text{ if } (x = x' \ \& \ y \ne y') \text{ or } (x \sim x' \ \& \ y \not\sim y')
\end{align*}
\end{definition}

In~\cite{Roberson12} it was further required that the $\rho_{xy}$ be projectors and that $\rho = I$. However, by replacing each $\rho_{xy}$ by the projection onto its image, it is easy to see these two definitions are equivalent. Note that in the above definition, $y \not\sim y'$ includes the case where $y = y'$, since we do not consider a vertex to be adjacent to itself.

\subsection{Entanglement-Assisted Homomorphisms}

Suppose that Alice wants to send Bob a message $x$ from a set $V$ of possible messages through the use of a noisy classical channel $\mathcal{N}$. Suppose further $\{C_1, C_2, \ldots C_m\}$ is a family of subsets of $V$ and that Bob is given some $i$ such that $x \in C_i$. Bob wants to recover Alice's message without error in every possible case of choices for $x$ and $i$. The noisy channel is specified by a conditional probability distribution, meaning that upon inputting an element $y$ of its input set, the channel outputs $w$ from its output set with probability $\mathcal{N}(w|y)$.

In order to analyze the above scenario, it is helpful to define two graphs. The \emph{disitinguishability graph of the channel}, which we will call $Y$, has the inputs of the channel as its vertices, and two vertices are adjacent if they cannot produce the same output. In other words, $y$ is adjacent to $y'$ if $\mathcal{N}(w|y)\mathcal{N}(w|y') = 0$ for all outputs $w$. The \emph{characteristic graph of Bob's side information}, which we will call $X$, has the message set $V$ as its vertex set, and two vertices $x$ and $x'$ are adjacent if there exists $i \in [m]$ such that $x,x' \in C_i$.

In~\cite{NTR} it was shown that classical Alice and Bob can perform the above communication task perfectly if and only if $X \to Y$. The idea is that messages which are not distinguishable by Bob using his side information, and thus are adjacent in $X$, must be mapped to channel inputs which are distinguishable, and thus adjacent in $Y$. Later it was shown that sharing an entangled state can sometimes allow Alice and Bob to perform perfectly even when $X \not\to Y$~\cite{ehomos}. However, it was also shown that the existence of a perfect quantum strategy still only depends on the graphs $X$ and $Y$. Thus we say that $X$ has an \emph{entanglement-assisted homomorphism} to $Y$ if there is a perfect quantum strategy for the above communication task. Mathematically we have the following~\cite{ehomos}:

\begin{definition}
We say that there exists a \emph{entanglement-assisted homomorphism} from a graph $X$ to a graph $Y$, and write $X\earrow Y$, if there exist positive semidefinite matrices $\rho_{xy}$ for $x \in V(X), y \in V(Y)$ and $\rho$ satisfying
\begin{align*}
&\langle \rho, \rho \rangle = 1 \\
&\sum_{y \in V(Y)} \rho_{xy} = \rho \text{ for all } x \in V(X) \\
&\langle \rho_{xy}, \rho_{x'y'} \rangle = 0 \text{ if } (x \sim x' \ \& \ y \not\sim y')
\end{align*}
\end{definition}

Clearly, the definition for entanglement-assisted homomorphism is weaker than that of quantum homomorphism, and so $X \qarrow Y$ implies $X \earrow Y$. Interestingly, it is not known whether or not the converse holds. What is known is that the converse is equivalent to projective measurements and maximally entangled states always sufficing for perfect quantum strategies to the above communication task. In terms of the mathematical definitions, the converse is equivalent to projectors always sufficing when choosing $\rho$ and the $\rho_{xy}$ in the definition of entanglement-assisted homomorphism.

\section{Cones}

We will specifically be concerned only with cones of matrices in this work, and so this is what we will implicitly mean when we say ``cone". Below we will give definitions of conic homomorphisms which could be applied to any cone, but if the cone does not have certain properties the corresponding homomorphism may be degenerate or not well-behaved. Because of this, here we will establish the properties we require of a cone to which we apply the below definitions.

One of the properties we will require of suitable cones is that of convexity. However, there is a subtlety here that needs to be pointed out. A cone is convex if it is closed under linear combinations with positive coefficients. For example, the set of $n \times n$ positive semidefinite matrices, denoted $\psdn$, is a convex cone. However, the set of all positive semidefinite matrices of any size, denoted $\psd$, is closed under multiplication by positive scalars and therefore still a cone, but cannot strictly be said to be convex since one cannot add matrices of different sizes. Still, since we will be dealing with matrices of various sizes, for our purposes we will refer to $\psd$ as being convex. Similarly, if $\K$ is the disjoint union of the convex cones $\K^n$, then we will say that $\K$ is convex.

One of the properties that we will need our cones to satisfy requires a little explanation. Given an $n \times n$ matrix $M$ and a partition $P = \{P_1, P_2, \ldots, P_m\}$ of $[n]$, we say that the matrix $N$ given by
\[N_{i,j} = \sum_{\ell \in P_i, k \in P_j} M_{\ell,k}\]
is the contraction of $M$ with respect to partition $P$. Contractions will be part of several of our proofs and thus we will require that our cones are closed under this operation.


We are now able to give the following definition.

\begin{definition}
A cone $\K$ is \emph{frabjous} if the following hold:
\begin{enumerate}
\item $\cp \subseteq \K \subseteq \psd$;\label{nested}
\item $\K$ is convex;\label{convex}
\item $\K$ is closed under taking principal submatrices;\label{submat}
\item $\K$ is closed under Kronecker product;\label{tensor}
\item $\K$ is closed under conjugation by permutation matrices;\label{perm}
\item $\K$ is closed under contraction.\label{contract}
\end{enumerate}
\end{definition}

This may seem like a lot of conditions to impose, but in practice they are not so unreasonable. In particular, the four cones we will be primarily interested in are frabjous cones. For some of our definitions we will want to further restrict to cones which contain only entrywise nonnegative matrices. For brevity we will refer to these as nonnegative cones.

\subsection{Cones of Interest}

There are four specific cones whose conic homomorphisms will correspond to previously defined notions, and for this reason we will pay them special attention. They are the following:

\begin{itemize}
\item $\cp$ - The \emph{completely positive} cone consists of the Gram matrices of entrywise nonnegative vectors.
\item $\cpsd$ - The \emph{completely positive semidefinite} cone consists of the Gram matrices of positive semidefinte matrices.
\item $\dnn$ - The \emph{doubly nonnegative cone} consists of the positive semidefinite matrices which are also entrywise nonnegative. Alternatively, the Gram matrices of vectors which have pairwise nonnegative inner products.
\item $\psd$ - The \emph{positive semidefinite} cone which consists of the Gram matrices of arbitrary vectors.
\end{itemize}

The Gram matrix formulations for these cones have been stressed because it will be relevant in some of our proofs. A Gram matrix of positive semidefinite matrices is likely a new concept for many, so we will clarify this definition. The Gram matrix of positive semidefinite matrices $\rho_1, \ldots, \rho_n$ is the matrix $M$ given by $M_{i,j} = \tr(\rho_i\rho_j)$.

It is not hard to see that the above cones are nested:
\[\cp \subset \cpsd \subset \dnn \subset \psd,\]
and it is known that all of these containments are strict~\cite{Piovesan14}.

\section{Conic Homomorphisms}

We define two different types of conic homomorphisms:

\begin{definition}
For graphs $X$ and $Y$, and nonnegative frabjous cone $\K$, we say that $X$ has a \emph{weak $\K$-homomorphism} to $Y$, and write $X \tokk Y$, if there exists a matrix $H \in \K$ satisfying
\begin{align}
\sum_{y,y'} H_{xy,x'y'} &= 1 \text{ for all } x,x' \in V(X) \label{cond:sum}\\
H_{xy,x'y'} &= 0 \text{ if } x \sim x' \text{ and } y \not\sim y' \label{cond:ortho}
\end{align}
Further, for any frabjous cone $\K$, we say that $X$ has a \emph{strong $\K$-homomorphism}, and write $X \tok Y$ if there exists $H \in \K$ satisfying the above conditions as well as
\begin{align}
H_{xy,x'y'} = 0 \text{ for } x = x' \text{ and } y \ne y' \label{cond:mortho}
\end{align}
\end{definition}

We will refer to a matrix $H$ fulfilling one of the definitions above as a \emph{homomorphism matrix}. It is useful to think of such a matrix as having blocks indexed by the vertices of $X$. Then the first condition above says that the entries of each of these blocks sum to 1.

We do not restrict the definition of weak $\K$-homomorphism to nonnegative frabjous cones without reason. When attempting to extend this definition to general frabjous cones, certain complications arise which prevent us from doing so in a way which allows our results to continue to hold. To see an in depth discussion of these complications see Section~\ref{sec:weakhomo}.

Obviously, $X \tok Y$ implies $X \tokk Y$ for any frabjous cone $\K$. Furthermore, if $\K \subseteq \K'$, then the existence of a strong/weak $\K$-homomorphism implies the existence of a strong/weak $\K'$-homomorphism. These are simple but fundamental observations that will be used implicitly as we proceed down the rabbit hole.

\subsection{Returning to $\cp$, $\cpsd$, $\dnn$, and $\psd$}

Before we investigate properties of conic homomorphisms in general, we will turn our attention to what results from applying the above definitions to the four specific cones we discussed above. To begin we will show that both strong and weak $\cp$-homomorphisms are equivalent to the usual graph homomorphism.

\begin{theorem}
For any graphs $X$ and $Y$, the following are equivalent:
\begin{enumerate}
\item $X \to Y$;
\item $X \tocpt Y$;
\item $X \tocpb Y$.
\end{enumerate}
\end{theorem}
\proof
Our proof will follow a $(1) \Rightarrow (2) \Rightarrow (3) \Rightarrow (1)$ structure.

$(1) \Rightarrow (2)$: Let $\varphi: V(X) \to V(Y)$ be a homomorphism from $X$ to $Y$. Define a vector $v \in \mathbb{R}^{V(X) \times V(Y)}$ as follows:
\[v_{xy} = 
\left\{\begin{array}{ll}
1 & \text{if } y = \varphi(x) \\
0 & o.w.
\end{array}\right.\]
It is easy to see that letting $H = vv^\intercal$ gives a strong $\cp$-homomorphism from $X$ to $Y$.

$(2) \Rightarrow (3)$: Trivial.

$(3) \Rightarrow (1)$: Suppose that $H$ is a homomorphism matrix certifying a weak $\cp$-homomorphism from $X$ to $Y$. Since $H \in \cp$, it is the Gram matrix of nonnegative vectors and thus it can be written as
\[H = P^\intercal P,\]
where $P$ is a nonnegative matrix. This simply says that $H$ is the Gram matrix of the columns of $P$. However, if we let $v_i$ be the $i^\text{th}$ row of $P$, then we have
\[H = P^\intercal P = \sum_i v_i v_i^\intercal.\]
For each $i$ let $H^i = v_iv_i^\intercal$. Each of these matrices is completely positive, and we aim to show that, up to a positive scalar, they are each a valid homomorphism matrix for $X \tocpb Y$. Since the $H^i$ are all nonnegative and $H$ satisfies condition~(\ref{cond:ortho}), the $H^i$ must all satisfy condition~(\ref{cond:ortho}). So we only need to check that $\sum_{y,y'} H^i_{xy,x'y'}$ is constant for all $x,x' \in V(X)$. This is equivalent to the contraction $M^i$, defined by
\[M^i_{x,x'} = \sum_{y,y'} H^i_{xy,x'y'}\]
being a scalar multiple of the all ones matrix $J$. Let $M$ be the contraction of $H$ defined analogously to the $M^i$. Clearly, $M = J$ and therefore
\[J = \sum_i M^i.\]
However, $J$ is a rank one positive semidefinite matrix and the $M^i$ are all positive semidefinite. Therefore each $M^i$ must be a scalar multiple of $J$.

From the above, we can assume that $H = vv^\intercal$ for some nonnegative vector $v \in \mathbb{R}^{V(X) \times V(Y)}$. Since
\[1 = \sum_{y,y'} H_{xy,xy'} = \sum_{y,y'} v_{xy} v_{xy'}\]
for all $x \in V(X)$, there exists some function $\varphi : V(X) \to V(Y)$ such that $v_{x\varphi(x)} \ne 0$ for all $x \in V(X)$. We claim $\varphi$ is a homomorphism. Indeed, if $x \sim x'$, then $H_{x\varphi(x),x'\varphi(x')} = v_{x\varphi(x)}v_{x'\varphi(x')} \ne 0$ implies that $\varphi(x) \sim \varphi(x')$. Therefore, if $X \tocpb Y$, then $X \to Y$ and this completes the proof.\qedn

%

Below we will see that the equivalence between $\tocpt$ and $\tocpb$ can be shown independently using techniques of Laurent and Piovesan from~\cite{Piovesan14}. Furthermore, the equivalence of $\tocpt$ and $\to$ will follow from a result of de Klerk and Pasechnik~\cite{Pasechnik02} and a result of ours that we prove in Section~\ref{sec:homotheta}. It is nice to have a direct proof of the equivalence of all three though.\\

In~\cite{Stahlke14}, three relaxations of entanglement-assisted homomorphisms were introduced. The first among these was denoted $\tob$ and was defined exactly as our definition of quantum homomorphism above, except that $\rho$ and the $\rho_{xy}$ were required to be vectors instead of positive semidefinite matrices. They also defined $\tov$ in the same manner with the added restriction that all inner products among the vectors must be nonnegative. Lastly, they defined $\toplus$ as we defined entanglement-assisted homomorphism above, except using vectors, and adding the same nonnegativity constraint on the inner products as for $\tov$. Note that this nonnegativity constraint would have been redundant in the definitions of quantum/entanglement-assisted homomorphism since the inner product of positive semidefinite matrices is always nonnegative.

It was also shown in~\cite{Stahlke14} that
\[X \tov Y \Rightarrow X \toplus Y \Rightarrow X \tob Y.\]
We will see later that $\tov$ and $\toplus$ are in fact equivalent. Interestingly, though they do it in different terms, Theorem 3 of~\cite{Stahlke14} proves that these three relations are equivalent to three of our conic homomorphisms:

\begin{remark}
The following pairs of relations are equivalent:
\begin{itemize}
\item $\tov$ and $\todnnt$;
\item $\toplus$ and $\todnnb$;
\item $\tob$ and $\topsdt$.
\end{itemize}
\end{remark}

It is worth noting that it is also shown in~\cite{Stahlke14} that $X \tob Y$ is equivalent to $\bartheta(X) \le \bartheta(Y)$, where $\bartheta$ is the Lov\'{a}sz theta number of the complement.

The main idea in showing the equivalences above is the following lemma which is implicitly proven in~\cite{Stahlke14}.

\begin{lemma}\label{lem:Gram}
Let $X$ and $Y$ be finite sets and let $H$ be the Gram matrix of vectors $w_{xy}$ for $x \in X$, $y \in Y$ which lie in some inner product space. Then there exists a vector $w$ such that 
\begin{align}\label{eqn:w}
\begin{split}
\langle w, w \rangle &= 1 \\ 
\sum_{y \in Y} w_{xy} &= w \text{ for all } x \in X
\end{split}
\end{align}
if and only if
\begin{equation}\label{eqn:sum}
\sum_{y,y' \in Y} H_{xy,x'y'} = 1 \text{ for all } x,x' \in X.
\end{equation}
\end{lemma}
\proof
For each $x \in X$, let $w_x = \sum_y w_{xy}$. Equation~(\ref{eqn:w}) is equivalent to the statement that the $w_x$ vectors are the same for all $x \in X$ and have norm 1. On the other hand, equation~(\ref{eqn:sum}) is equivalent to saying $\langle w_x, w_{x'} \rangle = 1$ for all $x,x' \in X$. The former clearly implies the latter. To see that converse, note that the latter says that the $w_x$ are all unit norm vectors with pairwise unit inner products. This is clearly only possible if they are equal.\qedn

Using the above lemma, the following theorem is immediate:

\begin{theorem}
The relations $\qarrow$ and $\tocpsdt$ are equivalent, as are $\earrow$ and $\tocpsdb$.
\end{theorem}
\proof
Matrices are vectors. Apply the above lemma.\qedn

So we have seen that the conic homomorphisms corresponding to the cones $\cp$, $\cpsd$, $\dnn$, and $\psd$ all correspond to previously defined notions. We also showed that $\tocpt$ and $\tocpb$ are equivalent, and now we will see that the same holds for $\dnn$.

\begin{theorem}\label{thm:dnnequiv}
Let $X$ and $Y$ be graphs. Then $X \todnnt Y$ if and only if $X \todnnb Y$.
\end{theorem}
\proof
Clearly, if $X \todnnt Y$ then $X \todnnb Y$. Conversely, suppose that $X \todnnb Y$ and that $H$ is the requisite homomorphism matrix. If $H_{xy,xy'} = 0$ for all $y \ne y'$ and $x \in V(X)$, then we are done. Suppose that $H_{xy,xy'} = c \ne 0$ for some $x \in V(X)$ and $y \ne y' \in V(Y)$. Then since $H \in \dnn$, we have that $c > 0$. Let $B$ be the matrix with $B_{xy,xy} = B_{xy',xy'} = 1 = -B_{xy,xy'} = -B_{xy',xy}$ and all other entries 0. Then $H + cB$ is the sum of two positive semidefinite matrices and is therefore positive semidefinite. Also, all entries of $H+cB$ are still nonnegative, therefore $H+cB \in \dnn$. Furthermore, we have removed one violation of condition~(\ref{cond:mortho}) and added none, and we have not caused any other conditions to be broken. Repeating this as necessary we obtain a matrix $H'$ which satisfies all of the requirements of $X \todnnt Y$.\qedn

The above theorem is a straightforward generalization of results from~\cite{Piovesan14} in which the above was done for the independence and chromatic number cases. They also use the same technique to show that $\tocpt$ and $\tocpb$ are equivalent in the special cases of independence and chromatic numbers. Importantly though, they need the following theorem from~\cite{cplemma}.

\begin{theorem}\label{thm:cp}
Assume that $A$ is completely positive, $B$ is positive semidefinite with all its entries equal to zero except for a $2 \times 2$ principal submatrix, and that $A + B$ is a
nonnegative matrix. Then A + B is completely positive.
\end{theorem}

If the above theorem held for completely positive semidefinite matrices as well, then we would be able to show that $\tocpsdt$ and $\tocpsdb$ are equivalent. However, the above theorem in fact does not hold for $\cpsd$, and a counterexample was given in~\cite{Piovesan14}.\\

From the equivalences we have shown above, we have reduced the number of distinct conic homomorphisms associated to the cones \ $\cone$ from seven to five. Furthermore, since $X \tocpsdb Y$ clearly implies $X \todnnb Y$ we can compare these five conic homomorphisms straightforwardly as follows:
\[X \tocpt Y \Rightarrow X \tocpsdt Y \Rightarrow X \tocpsdb Y \Rightarrow X \todnnt Y \Rightarrow X \topsdt Y.\]
From results of~\cite{Stahlke14} and others, the only two implications above whose converses may hold are $X \tocpsdt Y \Rightarrow X \tocpsdb Y \Rightarrow X \tov Y$, but it is not possible for the converses of both of these to hold.

In general, using results from~\cite{Stahlke14}, for any frabjous cone $\K$ we have that
\[X \to Y \Rightarrow X \tok Y \Rightarrow X \topsdt Y \Leftrightarrow \bartheta(X) \le \bartheta(Y),\]
and for any nonnegative frabjous cone we have
\[X \to Y \Rightarrow X \tokk Y \Rightarrow X \todnnt Y \Rightarrow \begin{array}{c} \barthetam(X) \le \barthetam(Y)  \\ \bartheta(X) \le \bartheta(Y) \\ \barthetap(X) \le \barthetap(Y)\end{array}\]
where $\barthetam$ and $\barthetap$ are Schrijver and Szegedy's variants of Lov\'{a}sz theta respectively.

\subsection{Basic Properties of Conic Homorphisms}

There are a few basic properties of conic homomorphisms that will be important throughout this paper, and that correspond to properties of usual graph homomorphisms. The first such property is that of reflexivity:

\begin{theorem}\label{thm:reflex}
For any frabjous cone $\K$ and graph $X$, we have that $X \tok X$. If $\K$ is also nonnegative, then also $X \tokk X$.
\end{theorem}
\proof
Define $\Phi = \sum_{x \in V(X)} e_x \otimes e_x$ where $e_x \in \mathbb{R}^{V(X)}$ is the $x^\text{th}$ standard basis vector. The matrix $\Phi \Phi^\intercal \in \cp \subseteq \K$ suffices for the conic homomorphism in both cases.\qedn

The vector $\Phi$ in the above proof is known as the maximally entangled state and is of considerable interest in the field of quantum information.

The next property is that of transitivity, though this requires a more involved proof.

\begin{theorem}\label{thm:trans}
Let $\K$ be a frabjous cone and let $X$, $Y$, and $Z$ be graphs. If $X \tok Y$ and $Y \tok Z$, then $X \tok Z$. If $\K$ is nonnegative as well, then the analogous statement holds for $\tokk$.
\end{theorem}
\proof
We will give the proof for $\tok$, but the proof for $\tokk$ is essentially the same. Let $H^1$ and $H^2$ be the homomorphism matrices certifying $X \tok Y$ and $Y \tok Z$ respectively. Define the matrix $H$ entrywise as follows:
\[H_{xz,x'z'} = \sum_{y,y'} H^1_{xy,x'y'} H^2_{yz,y'z'} = \sum_{y,y'} (H^1 \otimes H^2)_{xyyz,x'y'y'z'}.\]
The matrix $H$ is a contraction of a principal submatrix of $H^1 \otimes H^2$ and thus $H \in \K$. So we only need to check conditions~(\ref{cond:sum}), (\ref{cond:ortho}), and~(\ref{cond:mortho}).

We have that
\begin{align*}
\sum_{z,z'} H_{xz,x'z'} &= \sum_{z,z'} \sum_{y,y'} H^1_{xy,x'y'} H^2_{yz,y'z'} \\
&= \sum_{y,y'} H^1_{xy,x'y'} \sum_{z,z'} H^2_{yz,y'z'} \\
&= \sum_{y,y'} H^1_{xy,x'y'} = 1
\end{align*}
for all $x,x' \in V(X)$.

Now suppose that $x \sim x'$ and $z \not\sim z'$. For $y \not\sim y'$, we have that $H^1_{xy,x'y'}  = 0$, and for $y \sim y'$ we have that $H^2_{yz,y'z'} = 0$. Therefore
\[H_{xz,x'z'} = \sum_{y,y'} H^1_{xy,x'y'} H^2_{yz,y'z'} = 0.\]
Similarly, for all $x \in V(X)$ and $z \ne z' \in V(Z)$, we have that
\[H_{xz,xz'} = \sum_{y,y'} H^1_{xy,xy'} H^2_{yz,y'z'} = \sum_{y} H^1_{xy,xy} H^2_{yz,yz'} = 0.\]
Therefore $X \tok Z$ as desired.\qedn

If we let $a = |V(X)|$, $b = |V(Y)|$, and $c = |V(Z)|$, then the matrix $H$ constructed in the proof above can be written in a non-entrywise manner as follows:
\[H = (I_a \otimes \Phi_b^\intercal \otimes I_c)(H^1 \otimes H^2) (I_a \otimes \Phi_b \otimes I_c),\]
where $\Phi_b = \sum_1^b e_i \otimes e_i$. Perhaps this is enlightening in some way, but it does not seem to aid in the proof.

If the construction of the matrix $H$ in the above proof seems mysterious, keep in mind that it is actually motivated by the case where $\K = \cpsd$. Here, the matrices $H^1$ and $H^2$ are Gram matrices of some positive semidefinite matrices $\rho^1_{xy}$ for $x \in V(X)$, $y \in V(Y)$ and $\rho^2_{yz}$ for $y \in V(Y)$, $z \in V(Z)$ respectively. Then the matrix $H$ is the Gram matrix of
\[\rho_{xz} = \sum_{y} \rho^1_{xy} \otimes \rho^2_{yz}.\]
This construction is, in the author's opinion, easier to understand, and once one has this they can see how to write $H$ in terms of $H^1$ and $H^2$ as in the above proof. It is then just a matter of checking that the proof works more generally. It may be helpful to the reader to keep in mind that this generalization from the Gram matrix case is how most of the constructions of matrices in this paper were derived.

\section{New Thetas}

In~\cite{chif} and~\cite{Piovesan14} they define variants on the Lov\'{a}sz theta function by taking two formulations for $\vartheta$ and replacing the condition that the matrix $M$ is positive semidefinite with the condition that it is in a given cone $\K$ of matrices. We give the definition of the two variants below:

\begin{equation*}\label{eqn:thetak}
\vartheta^\K(X) =  \begin{array}[t]{ll}
\sup & \langle M,J \rangle \\
\text{s.t.} & M_{x,x'} = 0 \text{ if } x \sim x' \\
 & \tr(M) = 1 \\
 & M \in \K
\end{array}
\quad
\Theta^\K(X) =  \begin{array}[t]{ll}
\inf & t \\
\text{s.t.} & N_{x,x} = t \ \forall x \in V(X)\\
 & N_{x,x'} = 0 \text{ if } x \sim x' \\
 & N - J \succeq 0 \\
 & N \in \K
\end{array}
\end{equation*}

Note that in the above programs we use sup and inf instead of max and min. This is quite intentional since the cone $\K$ may not be closed. Indeed, it is an open question whether the cone $\cpsd$ is closed.

The value of the above two parameters is known for $\K \in \{\cp, \dnn, \psd\}$ and is as follows:

\begin{itemize}
\item $\psd$: $\vartheta^\psd$ is the usual Lov\'{a}sz $\vartheta$, and $\Theta^\psd  = \bartheta$;
\item $\dnn$: $\vartheta^\dnn$ is Schrijver's $\varthetam$, and $\Theta^\dnn$ is Szegedy's $\varthetap$ of the complement, denoted $\barthetap$;
\item $\cp$: $\vartheta^\cp$ is the usual independence number, denoted $\alpha$~\cite{Pasechnik02}, and $\Theta^\cp$ is the fractional chromatic number, denoted $\chi_f$~\cite{chif}.
\end{itemize}

It may be helpful to the reader to think of $\thetak$ as a type of independence number and $\Thetak$ as a type of chromatic number. From the above definitions it is easy to see that if $\K$ and $\K'$ are frabjous cones such that $\K \subseteq \K'$, then
\[\thetak(X) \le \vartheta^{\K'}(X) \text{ and } \Theta^{\K'}(X) \le \Thetak(X)\]
for any graph $X$. If we let $\barthetak(X) := \thetak(\overline{X})$ where $\overline{X}$ denotes the complement of $X$, then we further have that
\[\omega(X) \le \barthetak(X) \le \bartheta(X) \le \Thetak(X) \le \chi_f(X).\]

The parameters $\thetacpsd$ and $\Thetacpsd$ were introduced in~\cite{Piovesan14} and are perhaps genuinely new graph parameters. However, in Section~\ref{sec:projpack} we will see another pair of parameters that appear closely related to these and may even be equal to them.

Interestingly, the three pairs of parameters listed above share a common relation. In particular, if $\K \in \conesm$ and $X$ is a graph, then it is known that
\[\thetak(X) \Thetak(X) \ge |V(X)|\]
with equality if $X$ is vertex transitive. It is not surprising then that this holds also for $\K = \cpsd$, and indeed for any frabjous cone. The following theorem was originally proven in~\cite{chif}, but we will give a proof here for completeness and because the ideas in the proof will be used later.

\begin{theorem}\label{thm:thetas}
Let $\K$ be a frabjous cone, and $X$ be a graph. Then
\[\thetak(X)\Thetak(X) \ge |V(X)|\]
with equality if $X$ is vertex transitive.
\end{theorem}
\proof

Let $N$ be a solution to the program for $\Thetak(X)$ with value $t$ and let $n = |V(X)|$. Then $\tr(N) = tn$. Define $M = N/tn$. Clearly, $M$ is a feasible solution for $\thetak(X)$. Furthermore, since $N - J \succeq 0$, we have that
\[0 \le \frac{1}{tn}\langle N-J, J \rangle = \langle M,J \rangle - \frac{n^2}{tn} = \langle M,J \rangle - \frac{n}{t}.\]
Therefore, $M$ is a feasible solution for $\thetak(X)$ with value at least $n/t$. Therefore $\thetak(X)\Thetak(X) \ge |V(X)|$.

Now suppose that $X$ is vertex transitive. Any permutation $\pi$ of $V(X)$ can be associated with a permutation matrix $P_\pi$. If $A$ is the adjacency matrix of $X$, then $\pi$ is an automorphism of $X$ if and only if
\[P_\pi^\intercal A P_\pi = A.\]
Let $M$ be a feasible solution for $\thetak(X)$. Since $M$ has a zero entry everywhere $A$ has a one, and since $\K$ is closed under conjugation by permutation matrices, we have that $P_\pi^\intercal M P_\pi$ is also a feasible solution for any automorphism $\pi$ of $X$. Furthermore, this solution has the same objective value as the original. Using this we can use $M$ to construct the following solution:
\[\overline{M} = \frac{1}{|\aut(X)|}\sum_{\pi \in \aut(X)} P_\pi^\intercal M P_\pi,\]
where $\aut(X)$ is the automorphism group of $X$.
Since $X$ is vertex transitive, the matrix $\overline{M}$ has constant diagonal and constant row sum. Let $t = \langle \overline{M}, J \rangle$. Then $\overline{M}_{x,x} = 1/n$ for all $x \in V(X)$ and each row sum of $\overline{M}$ is equal to $t/n$. Let $N = (n^2/t)\overline{M}$. We will show that $N$ is a feasible solution to $\Thetak(X)$ with the appropriate value. Clearly, $N_{x,x'} = 0$ when $x \sim x'$, and $N_{x,x} = n/t$. All that is left to show is that $N - J \succeq 0$. To see this note that $N$ has constant row sum, and therefore its maximum eigenvalue has the all ones vector as an eigenvector and all other eigenvectors are orthogonal to the all ones vector. This means that all of the eigenvectors of $N$ are eigenvectors of $J$ and so we can easily determine the eigenvalues of $N-J$. The row sum of $N$ is $n$ and thus the eigenvalue of $N - J$ associated to the all ones vector is $n - n = 0$. All other eigenvalues of $J$ are 0 and thus the remaining eigenvalues of $N-J$ are the other eigenvalues of $N$ which are nonnegative since $N \in \K \subseteq \psd$. Therefore $N-J \succeq 0$ and $N$ is a feasible solution to $\Thetak(X)$ with value $n/t$. This implies that $\thetak(X)\Thetak(X) \le |V(X)|$, and with the other inequality above we obtain the equality.\qedn

The above proof actually only requires a few properties which are possessed by frabjous cones, and so it in fact holds more generally. It was indeed proved in this greater generality by the authors of~\cite{chif}.

{\textbf{Note:}} The only thing we needed to show equality in the above proof is that we can always find a solution to $\thetak(X)$ which has constant diagonal and constant row sum. We will need this later in the proof of Theorem~\ref{thm:hom2theta}.

\section{Relating the Homomorphisms and Thetas}\label{sec:homotheta}

Here we will prove the analog of the classical result that $X \to Y$ if and only if $\alpha(X \ltimes Y) = |V(X)|$. The graph $X \ltimes Y$ is known as the homomorphic product of $X$ and $Y$. This graph has vertex set $V(X) \times V(Y)$ and $(x,y) \sim (x',y')$ if ($x = x'$ and $y \ne y$) or ($x \sim x'$ and $y \not\sim y'$). In other words, the adjacencies of $X \ltimes Y$ exactly correspond to the required 0 entries of the homomorphism matrix in the definition of a strong conic homomorphism from $X$ to $Y$. We will also prove in this section some monotonicity theorems for $\thetak$ and $\Thetak$ with respect to $\tok$ and $\tokk$.

\begin{theorem}\label{thm:hom2theta}
Let $X$ and $Y$ be graphs and $\K$ a frabjous cone. Then $X \tok Y$ if and only if $\thetak(X \ltimes Y) = |V(X)|$ and this value is attained by some feasible solution. If $Y$ is vertex transitive, then both statements are also equivalent to $\Thetak(X \ltimes Y) = |V(Y)|$ and this value being attained by some solution.
\end{theorem}
\proof
First note that since the vertices $\{(x,y) : y \in V(Y)\}$ form a clique for each $x \in V(X)$, we have that $\chi(\overline{X \ltimes Y}) \le |V(X)|$. Since $\thetak(Z) \le \vartheta(Z) \le \chi(\overline{Z})$ for any graph $Z$ and frabjous cone $\K$, we have that $\thetak(X \ltimes Y) = |V(X)|$ and this value being attained is equivalent to the existence of a solution with this value. Similarly, we have that $\omega(X \ltimes Y) \ge |V(Y)|$ and $\omega(Z) \le \Thetak(Z)$ for any graph $Z$ and frabjous cone $\K$. Thus $\Thetak(X \ltimes Y) = |V(Y)|$ and this value being attained is equivalent to the existence of solution of this value.

Suppose that $X \tok Y$ and that $H$ is the corresponding homomorphism matrix. By conditions~(\ref{cond:sum}) and~(\ref{cond:mortho}), we have that $\tr(H) = |V(X)|$ and that $\langle H, J \rangle = |V(X)|^2$. Furthermore, $H_{xy,x'y'} = 0$ whenever $(x,y) \sim (x',y')$ in $X \ltimes Y$. Therefore, $H/|V(X)|$ is a solution to $\thetak(X \ltimes Y)$ of value $|V(X)|$.

Conversely suppose that $M$ is a solution to $\thetak(X \ltimes Y)$ of value $|V(X)|$. Clearly $M$ satisfies conditions~(\ref{cond:ortho}) and~(\ref{cond:mortho}). So it suffices to show that $\sum_{y,y'} M_{xy,x'y'}$ is constant for all $x,x' \in V(X)$. If this is the case then by scaling $M$ we will obtain a valid homomorphism matrix for $X \tok Y$.

Let $n = |V(X)|$, and define the matrix $\widehat{M}$ entrywise as
\[\widehat{M}_{x,x'} = \sum_{y,y'} M_{xy,x'y'}.\]
So we must show that $\widehat{M}$ is a scalar multiple of the all ones matrix. Note that $\widehat{M}$ is a contraction of $M$ and is therefore positive semidefinite. Also, since $(x,y) \sim (x,y')$ for $y \ne y'$, we have that
\[\tr(\widehat{M}) = \tr(M) = 1.\]
Furthermore, the sum of the entries of $\widehat{M}$ is equal to the sum of the entries of $M$. Therefore, if $e$ is the constant vector of unit norm, we have that
\[e^\intercal \widehat{M} e = \frac{1}{n}\langle M, J \rangle = 1.\]
This implies that the maximum eigenvalue is at least 1. But from the trace we see that the sum of the eigenvalues is 1. Therefore $\widehat{M}$ has only one nonzero eigenvalue which is 1 and the eigenvector for this is the all ones vector. Therefore $\widehat{M}$ is, up to a scalar, the all ones matrix. Therefore, up to a positive scalar, $M$ is a valid homomorphism matrix for $X \tok Y$.

To prove the second claim, suppose that $\Thetak(X \ltimes Y) = |V(Y)|$ and this value is attained by some solution. By applying the transformation used in the proof of Theorem~\ref{thm:thetas} we can obtain a solution for $\thetak(X \ltimes Y)$ of value
\[\frac{|V(X \ltimes Y)|}{|V(Y)|} = |V(X)|\]
and thus $X \tok Y$ by the above argument. Note that this direction does not require the vertex transitivity of $Y$.

Suppose that $Y$ is vertex transitive and $X \tok Y$. By the above argument we can find a matrix $M$ which is a solution for $\thetak(X \ltimes Y)$ of value $|V(X)|$. We will show that we can assume that $M$ has constant diagonal and row sum, and therefore we can apply the same transformation as in the vertex transitive case of Theorem~\ref{thm:thetas} to obtain a solution for $\Thetak(X \ltimes Y)$ of value$|V(X \ltimes Y)|/|V(X)| = |V(Y)|$.

By using automorphisms of $Y$ as in the proof of Theorem~\ref{thm:thetas}, we can assume that $M_{xy,xy} = M_{xy',xy'}$ for all $x \in V(X)$, $y,y' \in V(Y)$. This implies that
\[M_{xy,xy} = \frac{1}{|V(X)|\cdot|V(Y)|}\]
for all $x$ and $y$, and thus $M$ has constant diagonal.

Since $M \in \K \subseteq \psd$, there exist vectors $w_{xy}$ for $x \in V(X)$, $y \in V(Y)$ such that $M_{xy,x'y'} = \langle w_{xy}, w_{x'y'} \rangle$. Furthermore, from the argument above we know that the sums
\[\sum_{y,y'} M_{xy,x'y'}\]
are constant for all $x,x' \in V(X)$. By Lemma~\ref{lem:Gram}, we have that there exists a vector $w$ such that
\[\sum_{y} w_{xy} = w\]
for all $x \in V(X)$. This implies that
\[\sum_{y'} M_{xy,x'y'} = \sum_{y'} \langle w_{xy}, w_{x'y'} \rangle = \langle w_{xy}, w \rangle,\]
which does not depend on $x'$. However, since $M_{xy,xy'} = 0$ for $y \ne y'$, we know that the above sum is equal to $M_{xy,xy} = 1/\left(|V(X)|\cdot|V(Y)|\right)$. Therefore,
\[\sum_{x',y'} M_{xy,x'y'} = |V(X)| \frac{1}{|V(X)|\cdot|V(Y)|} = |V(Y)|,\]
and thus $M$ has constant row sum. Since $M$ has constant diagonal and row sum, by the argument of Theorem~\ref{thm:thetas}, we can scale $M$ to obtain a solution for $\Thetak(X \ltimes Y)$ of value $|V(Y)|$.\qedn


Letting $\K = \cp$ in the above theorem we see that $X \tocpt Y$ if and only if $\alpha(X \ltimes Y) = |V(X)|$. If we had not already done so this would have shown that $\tocpt$ is equivalent to $\to$, since it is widely known that the latter is equivalent to $\alpha(X \ltimes Y) = |V(X)|$. For $\K \in \{\dnn, \psd\}$, the above theorem was shown in~\cite{bacik95}. Finally, this result was proven for $\K = \cpsd$ recently by Stahlke~\cite{privatedan}. Even though the above was previously known for all $\K \in \cones$, it is enlightening to see that these are all essentially the same result.

\subsection{Monotonicity Theorems}

We say that a graph parameter $\beta$ is \emph{homomorphism monotone} if $X \to Y$ implies $\beta(X) \le \beta(Y)$ for any graphs $X$ and $Y$. Many homomorphism monotone graph parameters are known, such as clique, chromatic, and fractional chromatic numbers. The practical application of these results is often through the contrapositive, e.g.~if $\chi(X) > \chi(Y)$ then there exists no homomorphism from $X$ to $Y$. A similar type of application is that of obtaining bounds on other graph parameters. For instance, the bound $\bartheta(X) \le \chi(X)$ for all graphs $X$ follows from the monotonicity of $\bartheta$ and the fact that $\bartheta(K_n) = n$ for all $n$. More recently, results showing that $\barthetam$, $\bartheta$, and $\barthetap$ are all quantum and entanglement-assisted homomorphism monotone~\cite{Roberson12, Stahlke14} have allowed for bounds on quantum graph parameters such as quantum chromatic number. Here, through a general treatment of conic homomorphisms and the parameters $\thetak$ and $\Thetak$, we will recover very many of these monotonicity results, as well as prove some new ones.

\begin{theorem}\label{thm:monot}
Let $X$ and $Y$ be graphs and $\K$ a frabjous cone. If $X \tok Y$, then $\barthetak(X) \le \barthetak(Y)$.
\end{theorem}
\proof
First we point out that one should keep in mind that we are considering $\barthetak$, and therefore the required zero entries of the feasible solutions correspond to distinct nonadjacent vertices.

Suppose that $X \tok Y$ and $H$ is the corresponding homomorphism matrix. Let $M$ be a solution for $\barthetak(X)$. We will construct a feasible solution for $\barthetak(Y)$ of value equal to that of $M$. Define the matrix $N$ entrywise as follows:
\[N_{y,y'} = \sum_{x,x'} M_{x,x'} H_{xy,x'y'} = \sum_{x,x'} (M \otimes H)_{xxy,x'x'y'}.\]
Since $H,M \in \K$ and $N$ is a contraction of a principal submatrix of $M \otimes H$, we have $N \in \K$. Considering the trace of $N$ and first using the zero entries of $M$ and then of $H$, we see that
\begin{align*}
\tr(N) &= \sum_{x,x',y} M_{x,x'} H_{xy,x'y} \\
&= \sum_{x = x' \text{ or } x \sim x', \ y} M_{x,x'} H_{xy,x'y} \\
&= \sum_{x, y} M_{x,x} H_{xy,xy} \\
&= \sum_x M_{x,x} \sum_y H_{xy,xy} \\
&= \sum_x M_{x,x} = 1.
\end{align*}

Next, suppose that $y$ and $y'$ are distinct nonadjacent vertices of $Y$. Then, again first using the zero entries of $M$ and then of $H$, we have
\begin{align*}
N_{y,y'} &= \sum_{x,x'} M_{x,x'} H_{xy,x'y'} \\
&= \sum_{x = x' \text{ or } x \sim x'} M_{x,x'} H_{xy,x'y'} \\
&= 0.
\end{align*}
This implies that $N$ is indeed a feasible solution for $\barthetak(Y)$.

Lastly, the value of the solution $N$ is given by
\begin{align*}
\langle N, J \rangle &= \sum_{y,y'} N_{y,y'} = \sum_{x,x',y,y'} M_{x,x'} H_{xy,x'y'} \\
&= \sum_{x,x'} M_{x,x'} \sum_{y,y'} H_{xy,x'y'} \\
&= \sum_{x,x'} M_{x,x'} = \langle M, J \rangle.
\end{align*}\qedn

%

After the above theorem, it is natural to ask if something similar can be shown for $\Thetak$. It turns out we can prove an even stronger result for these parameters.

\begin{theorem}\label{thm:monob}
Let $X$ and $Y$ be graphs and $\K$ a frabjous cone. If $X \tok Y$, then $\Thetak(X) \le \Thetak(Y)$. Furthermore, if $\K$ is a nonnegative frabjous cone, then $X \tokk Y$ implies $\Thetak(X) \le \Thetak(Y)$.
\end{theorem}
\proof
The proofs of the two statements are similar, so we will just prove the latter. Suppose that $N$ is a solution for $\Thetak(Y)$ of value $N_{yy} = t$, and that $H$ is a homomorphism matrix for $X \tokk Y$. Define the matrix $M$ entrywise as follows:
\[M_{x,x'} = \sum_{y,y'} H_{xy,x'y'} N_{y,y'} = \sum_{y,y'} (H \otimes N)_{xyy,x'y'y'}.\]
Then $M$ is a contraction of a principal submatrix of $H \otimes N$, and therefore $M \in \K$. Since $N$ is positive semidefinite and has constant diagonal $t$, we have that $N_{y,y'} \le t$ for all $y,y' \in V(Y)$. Therefore, for $x \in V(X)$, we have
\[M_{x,x} = \sum_{y,y'} H_{xy,xy'} N_{y,y'} \le t \sum_{y,y'} H_{xy,xy'} = t.\]
Note that the above is the only place the proof would differ for the $X \tok Y$ case for a possibly not nonnegative frabjous cone. Since in that case $N_{y,y'}$ could possibly be negative, we would need to use the fact that condition~(\ref{cond:mortho}) implies that $H_{xy,xy'} \ge 0$ for all $y,y' \in V(Y)$, rather than simply using the fact that $\K$ is nonnegative.

So the diagonal entries of $M$ are all at most $t$, but we would like them to all be equal to $t$. In order to do this, we can simply add an appropriate nonnegative diagonal matrix. Since any such matrix is in $\cp \subseteq \K$, this will not move $M$ outside of $\K$. Furthermore, adding a nonnegative diagonal matrix to $M$ will not disrupt any of the other feasibility conditions for $\Thetak(X)$.

Now suppose that $x \sim x'$. Then, first using the zeros of $H$ and then of $N$, we have that
\[M_{x,x'} = \sum_{y,y'} H_{xy,x'y'} N_{y,y'} = \sum_{y \sim y'} H_{xy,x'y'} N_{y,y'} = 0.\]

All that is left is to show that $M - J$ is positive semidefinite. We have that
\begin{align*}
M_{x,x'} - 1 &= \left(\sum_{y,y'} H_{xy,x'y'} N_{y,y'}\right) - 1 \\
&= \left(\sum_{y,y'} H_{xy,x'y'} N_{y,y'}\right) - \left(\sum_{y,y'} H_{xy,x'y'}\right) \\
&= \sum_{y,y'} \left(H \otimes N\right)_{xyy,x'y'y'} - \sum_{y,y'} \left(H \otimes J\right)_{xyy,x'y'y'} \\
& = \sum_{y,y'} \left(H \otimes \left(N - J\right)\right)_{xyy,x'y'y'}
\end{align*}
So $M - J$ is a contraction of a principal submatrix of $H \otimes (N-J)$, and both $H$ and $N-J$ are positive semidefinite. Therefore $M - J$ is positive semidefinite and thus $M$ is a feasible solution for $\Thetak(X)$ of value $t$.\qedn

Of course for $\K \in \conesm$ the hypothesis here is the same as in the previous theorem. But it is interesting that $\Theta^\cpsd$ is monotone with respect to $\earrow$, but $\bartheta^\cpsd$ perhaps is not. If you try to prove that it is monotone using the technique above, you will see that you are not able to get the necessary zero entries.

A consequence of the above two theorems is that if $\K$ is any frabjous cone, then $X \tok Y$ implies that $\bartheta^{\K'}(X) \le \bartheta^{\K'}(Y)$ and $\Theta^{\K'}(X) \le \Theta^{\K'}(Y)$ for any frabjous cone $\K'$ satisfying $\K \subseteq \K'$. Similarly, if $\K$ is a nonnegative frabjous cone and $X \tokk Y$, then $\Theta^{\K'}(X) \le \Theta^{\K'}(Y)$ for any $\K'$ containing $\K$. Varying $\K$ and $\K'$, one recovers virtually all known monotonicity theorems regarding $\barthetam$, $\bartheta$, and $\barthetap$. Furthermore, we obtain the new results $X \qarrow Y \Rightarrow \bartheta^\cpsd(X) \le \bartheta^\cpsd(Y)$ and $X \earrow Y \Rightarrow \Theta^\cpsd(X) \le \Theta^\cpsd(Y)$.


\section{Projective packings and projective rank}\label{sec:projpack}

A projective packing is an assignment $x \mapsto P_x$ of Hermitian projectors to the vertices of a graph $X$ such that adjacent vertices receive orthogonal projectors. The value of a projective packing is equal to the sum of the ranks of the projectors divided by their dimension. Similarly, a $d/r$-projective representation of a graph $X$ is a projective packing using only rank $r$ projectors in dimension $d$, and $\frac{d}{r}$ is the value of the representation. The supremum of the values of projective packings of a graph $X$ is the \emph{projective packing number} of $X$ and is denoted $\alpha_p(X)$. The projective rank of $X$ is the infimum of $\frac{d}{r}$ such that $X$ admits a $d/r$-representation, and this is denoted $\xi_f(X)$. The projective packing number and projective rank are relaxations of the independence number and fractional chromatic number respectively. To see this for the projective packing number, consider the projective packing which assigns the identity to every vertex in a maximum independent set of $X$ and assigns the zero matrix everywhere else. This has a value of $\alpha(X)$ and thus $\alpha(X) \le \alpha_p(X)$. Similarly, one can show that $\xi_f(X) \le \chi_f(X)$~\cite{RobersonThesis}. These parameters were originally introduced in connection with quantum homomorphisms and it has been shown that $\alpha_p$ of the complement and $\xi_f$ are quantum homomorphism monotone~\cite{Roberson12, RobersonThesis}. Moreover, an analogous result to Theorem~\ref{thm:hom2theta} has been shown for $\alpha_p$ and $\xi_f$~\cite{RobersonThesis}.

We will see that the pair of parameters $\alpha_p$ and $\xi_f$ is very similar to the pair $\thetacpsd$ and $\Thetacpsd$. For instance, it is not hard to see that $\alpha_p(X)\xi_f(X) \ge |V(X)|$ with equality if $X$ is vertex transitive~\cite{RobersonThesis}. We also have the following relation:

\begin{theorem}
For any graph $X$,
\[\alpha_p(X) \le \vartheta^\cpsd(X).\]
\end{theorem}
\proof
Let $x \mapsto P_x \in \mathbb{R}^{d \times d}$ be a projective packing of $X$ and let $P = \sum_x P_x$. Since $\rk(P_x) = \tr(P_x)$, the value of the projective packing is $\tr(P)/d$. Let $M'$ be the Gram matrix of the $P_x$. Since $P_x^2 = P_x$ for all $x \in V(X)$, we have that $\tr(M') = \tr(P)$, and so letting $M = M'/\tr(P)$ gives a feasible solution for $\vartheta^\cpsd(X)$. The value of this solution is $\tr(P^2)/\tr(P)$ and thus we only need to show that this is at least $\tr(P)/d$. This is equivalent to
\[d\tr(P^2) \ge \tr(P)^2,\]
which follows by applying the Cauchy-Schwarz inequality to the vector of eigenvalues of $P$ and the all ones vector.\qedn

As one might expect, there is a similar relationship between $\xi_f$ and $\Theta^\cpsd$:

\begin{theorem}
For any graph $X$,
\[\Theta^\cpsd(X) \le \xi_f(X).\]
\end{theorem}
\proof
Let $x \mapsto P_x$ be a $d/r$-projective representation of $X$ and let $N$ be the Gram matrix of the $\frac{\sqrt{d}}{r}P_x$. Clearly, $N \in \cpsd$ and has the necessary zero entries and constant diagonal equal to $\frac{d}{r}$. All that is left check is that $N-J$ is positive semidefinite. It is easy to check that $N-J$ is the Gram matrix of the vectors
\[u_x = \frac{\sqrt{d}}{r}\vect(P_x) - \frac{1}{\sqrt{d}}\vect(I)\]
thus finishing the proof.\qedn

It is possible that the inequalities in the above two theorems actually hold with equality. For the former, if the matrix $P$ in the proof is not a multiple of identity, then the value of $M$ is strictly greater than that of the projective packing. So if $\alpha_p(X) = \vartheta^\cpsd(X)$ does hold, then every optimal projective packing (if any exist) must sum up to a multiple of identity. For the latter pair of parameters, it is interesting to note that Theorem~\ref{thm:monob} implies that $\Theta^\cpsd(X)$ is monotone with respect to entanglement-assisted homomorphisms, and yet the author does not see any way to show that this is the case for $\xi_f(X)$.

\section{Graph products}

The Lov\'{a}sz theta number was originally introduced as an upper bound on the Shannon capacity of a graph~\cite{lovasz}. This result followed from two properties of $\vartheta$: 1) it upper bounds the independence number, and 2) it is multiplicative with respect to the strong graph product (see definition below). It has since been shown that $\vartheta$ is multiplicative with respect to the disjunctive product~\cite{knuth}. These two results also imply that $\vartheta$ is multiplicative with respect to the lexicographic product, since this is nested between the strong and disjunctive products in the subgraph order.

Here we will show that $\thetak$ shares some of these multiplicativity properties, namely that it is multiplicative with respect to the disjunctive and lexicographic products, but not the strong product in general. For $\Thetak$, we are only able to show that it is submultiplicative with respect to the disjunctive product, but we suspect it is in fact multiplicative with respect to this and the lexicographic product. We will also investigate the categorical product in this section, showing that its use as a meet operation for the homomorphism order extends to conic homomorphisms.

\subsection{Multiplicativity}

For graphs $X$ and $Y$, the strong, lexicographic, and disjunctive products are denoted $X \boxtimes Y$, $X[Y]$, and $X * Y$ respectively. All of these products have $V(X) \times V(Y)$ as their vertex set, and vertices $(x,y)$ and $(x',y')$ are adjacent in each of them as follows:
\begin{itemize}
\item $X \boxtimes Y$: ($x \sim x'$ \& $y \sim y'$) or ($x \sim x'$ \& $y = y'$) or ($x = x'$ \& $y \sim y'$);
\item $X[Y]$: $x \sim x'$ or ($x = x'$ \& $y \sim y'$);
\item $X * Y$: $x \sim x'$ or $y \sim y'$.
\end{itemize}

We say that a graph parameter $\beta$ is \emph{multiplicative} with respect to the strong product if $\beta(X \boxtimes Y) = \beta(X)\beta(Y)$ for all graphs $X$ and $Y$. Multiplicativity with respect to other products is defined analogously. If the equality above is replaced with `$\le$' then we say that $\beta$ is \emph{submultiplicative}, and \emph{supermultiplicative} is similarly defined.

The proof of the following theorem showing that $\thetak$ is multiplicative with respect to the disjunctive and lexicographic products is essentially the same as the proof for $\varthetam$ given in~\cite{Stahlke14}. Here we simply note that it holds for any frabjous cone.

\begin{theorem}
Let $X$ and $Y$ be graphs and $\K$ a frabjous cone. Then
\[\thetak(X[Y])  = \thetak(X) \thetak(Y) = \thetak(X * Y).\]
\end{theorem}
\proof
The proof will proceed by showing that $\thetak(X[Y]) \le \thetak(X) \thetak(Y) \le \thetak(X * Y)$. Since $X[Y]$ is a subgraph of $X * Y$, we have that $\thetak(X * Y) \le \thetak(X[Y])$ and so equality must hold throughout in the previous inequalities.

We begin by showing that $\thetak(X * Y) \ge \thetak(X) \thetak(Y)$. Suppose that $M^1$ and $M^2$ are feasible solutions for $\thetak(X)$ and $\thetak(Y)$ respectively. We will show that $M = M^1 \otimes M^2$ is a feasible solution for $\thetak(X * Y)$. First, $\tr(M) = \tr(M^1 \otimes M^2) = 1 \cdot 1 = 1$. Next, if $(x,y) \sim (x',y')$, then either $x \sim x'$ or $y \sim y'$. In either case $M_{xy,x'y'} = M^1_{x,x'} M^2_{y,y'} = 0$. Therefore $M$ is feasible and its objective value is $\langle M, J \rangle = \langle M^1, J \rangle \langle M^2, J \rangle$. Therefore $\thetak(X * Y) \ge \thetak(X) \thetak(Y)$.

To see that $\thetak(X[Y]) \le \thetak(X)\thetak(Y)$, consider a solution $M$ for $\thetak(X[Y])$. For each $x \in V(X)$, define $M^x$ entrywise as
\[M^x_{y,y'} = M_{xy,xy'}.\]
In other words, $M^x$ is the diagonal block of $M$ corresponding to $x$. For $y \sim y'$, the vertices $(x,y)$ and $(x,y')$ are adjacent in $X[Y]$, and so $M^x/\tr(M^x)$ is a feasible solution for $\thetak(Y)$. This implies that
\[\frac{\langle M^x, J \rangle}{\tr(M^x)} \le \thetak(Y) \text{ for all } x \in V(X).\]
Also, since $\sum_x \tr(M^x) = \tr(M) = 1$, we have that
\[\sum_x \langle M^x, J \rangle \le \sum_x \thetak(Y) \tr(M^x) = \thetak(Y).\]
Now define the contraction $M^X$ of $M$ as follows:
\[M^X_{x,x'} = \sum_{y,y'} M_{xy,x'y'}.\]
Since $(x,y) \sim (x',y')$ whenever $x \sim x'$, the matrix $M^X/\tr(M^X)$ is a feasible solution for $\thetak(X)$. So we have that
\[\langle M^X, J \rangle \le \thetak(X) \tr(M^X).\]
Combining all of the above, we see that the objective value of $M$ is
\[\langle M, J \rangle = \langle M^X, J \rangle \le \thetak(X) \tr(M^X) = \thetak(X) \sum_x \langle M^x, J \rangle \le \thetak(X) \thetak(Y).\]
Therefore $\thetak(X[Y]) \le \thetak(X)\thetak(Y)$ and this completes the proof.\qedn

The above proof cannot be extended to the strong product. Indeed, in~\cite{Stahlke14} they give an example of graphs $X$ and $Y$ such that $\varthetam(X \boxtimes Y) > \varthetam(X)\varthetam(Y)$.

For $\Thetak$, we are only able to show that it is submultiplicative with respect to the disjunctive product.

\begin{theorem}
Let $X$ and $Y$ be graphs and $\K$ a frabjous cone. Then
\[\Thetak(X * Y) \le \Thetak(X) \Thetak(Y).\]
\end{theorem}
\proof
Let $M^1$ and $M^2$ be solutions for $\Thetak(X)$ and $\Thetak(Y)$ of objective values $t_1$ and $t_2$ respectively. Define $M = M^1 \otimes M^2$. It is easy to see that $M_{xy,x'y'} = 0$ whenever $(x,y) \sim (x',y')$, and that the objective value of $M$ is $t_1 t_2$. Therefore it is only left to show that $M - J \succeq 0$. To see this note that
\begin{align*}
M - J =& M^1 \otimes M^2 - J \otimes J \\
=& (M^1 \otimes M^2 - M^1 \otimes J - J \otimes M^2 + J \otimes J) \\
&+ (M^1 \otimes J - J \otimes J) + (J \otimes M^2 - J \otimes J) \\
=& (M^1 - J) \otimes (M^2 - J) + (M^1 - J) \otimes J + J \otimes (M^2 - J).
\end{align*}
Since $J \succeq 0$, and $M^1 - J, M^2 - J \succeq 0$ by assumption, the above shows that $M - J \succeq 0$.\qedn

Note that for $\K \in \{\cp, \psd\}$, it is known that $\Thetak$ is multiplicative to both the disjunctive and lexicographic product, and thus we expect that the same is true for all frabjous $\K$. On the other hand, it is known that $\Thetak$ cannot be multiplicative with respect to the strong product in general, since there are counterexamples to this for $\K = \dnn$ given in~\cite{Stahlke14}.

\subsection{Lattice Properties}

A major area in the field of graph homomorphisms is the study of the homomorphism order. This is a partial order on homomorphic equivalence classes of graphs induced by the relation $\to$. A key property of this partial order is that it is a lattice. The meet and join operations of this lattice are given by the categorical product and disjoint union of graphs respectively. The categorical product of graphs $X$ and $Y$, denoted $X \times Y$, is the graph on vertex set $V(X) \times V(Y)$ in which $(x,y) \sim (x',y')$ if $x \sim x'$ and $y \sim y'$. The property which allows this product to be used as a meet operation is the following: a graph $Z$ admits a homomorphism to each of the graphs $X$ and $Y$ if and only if $Z \to X \times Y$. Here we will show that this continues to hold for conic homomorphisms.

\begin{theorem}
Let $\K$ be a frabjous cone and let $X$, $Y$, and $Z$ be graphs. Then $Z \tok X \times Y$ if and only if $Z \tok X$ and $Z \tok Y$. If $\K$ is also nonnegative then the analogous statement holds for $\tokk$.
\end{theorem}
\proof
The proofs for strong and weak conic homomorphisms are sufficiently similar and so we will simply prove it for the former. Suppose that $Z \tok X \times Y$. The map $(x,y) \mapsto x$ is a homomorphism from $X \times Y$ to $X$, and therefore $X \times Y \to X$ which implies $X \times Y \tok X$. Then by transitivity we have that $Z \tok X$, and similarly $Z \tok Y$.

Now suppose that $H^1$ and $H^2$ respectively are homomorphism matrices for $Z \tok X$ and $Z \tok Y$. Define the matrix $H$ entrywise as follows:
\[H_{zxy,z'x'y'} = H^1_{zx,z'x'} H^2_{zy,z'y'} = (H^1 \otimes H^2)_{zxzy,z'x'z'y'}.\]
Since $H$ is a principal submatrix of $H^1 \otimes H^2$ we have that $H \in \K$. For any $z, z' \in V(Z)$ we have
\begin{align*}
\sum_{x,y,x',y'} H_{zxy,z'x'y'} &= \sum_{x,y,x',y'} H^1_{zx,z'x'} H^2_{zy,z'y'} \\
&= \left(\sum_{x,x'} H^1_{zx,z'x'}\right) \left(\sum_{y,y'} H^2_{zy,z'y'}\right) = 1.
\end{align*}
From the definition of $H$ it is straightforward to see that conditions~(\ref{cond:ortho}) and~(\ref{cond:mortho}) both hold for $H$.\qedn

Given the above theorem it is natural to ask whether the disjoint union of graphs can be used as a join operation for conic homomorphisms as it is for homomorphisms. The ability of the disjoint union to fullfill this role comes from the following property: graphs $X$ and $Y$ both admit a homomorphism to a graph $Z$ if and only if their disjoint union, denoted $X + Y$, admits a homomorphism to $Z$. This is trivial to see and thus one would expect that this property easily extends to conic homomorphisms. Quite surprisingly, we are unable to show that this is the case for all frabjous cones. However, we can prove it under certain conditions. In particular, we have the following theorem:

\begin{theorem}
Let $X$, $Y$, and $Z$ be graphs and $\K \in \cones$. Then $X \tok Z$ and $Y \tok Z$ if and only if $X + Y \tok Z$. If $\K \in \{\cp, \cpsd, \dnn\}$, then the analogous statement holds for $\tokk$.
\end{theorem}
Again, we will only prove the claim for strong conic homomorphisms, but the argument is the same for weak. If $X + Y \tok Z$, then since $X \to X + Y$ we have that $X \tok X + Y$ and thus $X \tok Z$ by transitivity, and similarly $Y \tok Z$.

Now suppose that $X \tok Z$ and $Y \tok Z$ and that $H^1$ and $H^2$ are the corresponding homomorphism matrices. Since $\K \in \cones$, we know that $H^1$ and $H^2$ can be written as the Gram matrices of some appropriate objects $\rho_{xz}$ for $x \in V(X)$, $z \in V(Z)$, and $\sigma_{yz}$ for $y \in V(Y)$ and $z \in V(Z)$. Furthermore, from condition~(\ref{cond:sum}) and Lemma~\ref{lem:Gram}, there exist correspondingly appropriate objects $\rho$ and $\sigma$ such that
\begin{align*}
\sum_{z} \rho_{xz} = \rho \text{ for all } x \in V(X); \\
\sum_{z} \sigma_{yz} = \sigma \text{ for all } y \in V(Y).
\end{align*}
For all $w \in V(X + Y)$ and $z \in V(Z)$, define
\[\tau_{wz} = 
\left\{\begin{array}{ll}
\rho_{wz} \otimes \sigma  & \text{if } w \in V(X) \\
\rho \otimes \sigma_{wz} & \text{if } w \in V(Y)
\end{array}\right.\]
Let $H$ be the Gram matrix of the $\tau_{wz}$. It is straightforward to check that for all four possibilities of $\K$, the construction given guarantees that $H \in \K$. If $w \in V(X)$ and $z \ne z' \in V(Z)$, then
\[H_{wz,wz'} = \langle \tau_{wz}, \tau_{wz'} \rangle = \langle \rho_{wz}, \rho_{wz'} \rangle \langle \sigma, \sigma \rangle = 0,\]
and similarly for $w \in V(Y)$. If $w \sim w'$ in $X + Y$, then either $w,w' \in V(X)$ or $w,w' \in V(Y)$. In the former case, if $z \not\sim z'$ then
\[H_{wz,w'z'} = \langle \tau_{wz}, \tau_{w'z'} \rangle = \langle \rho_{wz}, \rho_{w'z'} \rangle \langle \sigma, \sigma \rangle = 0,\]
and similarly for the latter case. Therefore $H$ satisfies conditions~(\ref{cond:ortho}) and~(\ref{cond:mortho}), and we only have left to check condtion~(\ref{cond:sum}).

For $w,w' \in V(X)$ we have that
\[\sum_{z,z'} H_{wz,w'z'} = \sum_{z,z'} \langle \tau_{wz}, \tau_{w'z'} \rangle = \sum_{z,z'} \langle \rho_{wz}, \rho_{w'z'} \rangle \langle \sigma, \sigma \rangle = \langle \rho, \rho \rangle \langle \sigma, \sigma \rangle = 1,\]
and similarly for $w,w' \in V(Y)$. Otherwise, without loss of generality we have that $w \in V(X)$ and $w' \in V(Y)$. In this case,
\[\sum_{z,z'} H_{wz,w'z'} = \sum_{z,z'} \langle \tau_{wz}, \tau_{w'z'} \rangle = \sum_{z,z'} \langle \rho_{wz}, \rho \rangle \langle \sigma, \sigma_{w'z'} \rangle = \langle \rho, \rho \rangle \langle \sigma, \sigma \rangle = 1.\]
Therefore $H$ satisfies condition~(\ref{cond:sum}) and thus $X + Y \tok Z$.\qedn

The matrix $H$ constructed in the above proof can be written purely in terms of $H^1$ and $H^2$ without resorting to their Gram representations. If we let $a = |V(X)|$, $b = |V(Y)|$, $c = |V(Z)|$, and $J$ be the $ac \times bc$ all ones matrix, then
\[H = \begin{pmatrix}
H^1 & \frac{1}{ab}H^1 J H^2 \\
\frac{1}{ab}(H^1 J H^2)^\dagger & H^2
\end{pmatrix}.\]

One way to extend the above result to all frabjous cones would be to change the definition of frabjous to directly require that the above matrix is in the cone whenever $H^1$ and $H^2$ are homomorphism matrices for conic homomorphisms to the same graph. This of course would be a very specific and seemingly unnatural condition to impose, which is why the author did not take this approach. Alternatively, one could attempt to show that this condition actually follows from the given conditions on frabjous cones. The author was unable to do this but we hold out hope that it is possible.

\section{Conic Parameters}

Graph homomorphisms can be used to define certain graph parameters such as clique and chromatic numbers, and we can similarly use conic homomorphisms to define conic analogs of these classical graph parameters. We will in particular be interested in the following:

\begin{definition}
For a frabjous cone $\K$, we define the strong and weak $\K$ independence numbers as follows:
\begin{align*}
\alpha^\K(X) &= \max\{n : K_n \tok \overline{X}\} \\
\alpha_\K(X) &= \max\{n : K_n \tokk \overline{X}\}
\end{align*}
where the latter is defined only for nonnegative cones.
\end{definition}

For $\K = \cp$, we of course just obtain the usual independence number. It was shown in~\cite{Stahlke14} that $\alpha^\psd = \lfloor \vartheta \rfloor$ and $\alpha^\dnn = \lfloor \varthetam \rfloor$. For $\K = \cpsd$, we obtain the quantum and entanglement-assisted independence numbers ($\alpha_q$ and $\alpha^*$) respectively. Note that by our monotonicity theorem for $\thetak$ we have that $\alpha^\K(X) \le \thetak(X)$ for any graph $X$, but the same may not be true for $\alpha_\K$.

In~\cite{Roberson12} it was shown that $X \qarrow Y$ if and only if $\alpha_q(X \ltimes Y) = |V(X)|$, or in terms of cones, $X \tocpsdt Y$ if and only if $\alpha^\cpsd(X \ltimes Y) = |V(X)|$. This is rather interesting in light of Theorem~\ref{thm:hom2theta}, since it is not known whether or not $\alpha^\K(X) = \lfloor \thetak(X) \rfloor$ for $\K = \cpsd$, as it does for $\K \in \conesm$. In this section we will show that the above relation between $\alpha^\cpsd$ and $\tocpsdt$ holds more generally for any frabjous cone $\K$. To do this, we will use the following curious lemma.

\begin{lemma}\label{lem:curious}
Let $X$ be a graph and $\K$ a frabjous cone. Then $\thetak(X) = \chi(\overline{X})$ and this value is attained by some solution if and only if $\alpha^\K(X) = \chi(\overline{X})$.
\end{lemma}
\proof
First, we have that $\alpha^\K(X) \le \thetak(X) \le \vartheta(X) \le \chi(\overline{X})$, and therefore if $\alpha^\K(X) = \chi(\overline{X})$ then $\thetak(X) = \chi(\overline{X})$. Furthermore, by the proof of Theorem~\ref{thm:monot} we see that a solution for $\thetak(X)$ of value $\alpha^\K(X)$ can always be obtained. This proves the backwards direction of the lemma.

Conversely, suppose that $M$ is a solution for $\thetak(X)$ of value $c := \chi(\overline{X})$, and let $\{S_i : i \in [c]\}$ be a set of $c$ cliques partitioning $V(X)$. Note that $M_{x,x'} = 0$ whenever $x,x' \in S_i$ and $x \ne x'$.  Therefore, if we define $N$ to be the $c \times c$ matrix defined by
\[N_{i,j} = \sum_{x \in S_i, x' \in S_j} M_{x,x'},\]
then $\tr(N) = \tr(M) = 1$ and $\langle N,J \rangle = \langle M,J \rangle = c$. So $N$ is a $c \times c$ positive semidefinite matrix with trace 1 and sum $c$. This implies that $N = J$ as in the proof of Theorem~\ref{thm:hom2theta}. From this we see that
\begin{equation}\label{eqn:M}
\sum_{x \in S_i, x' \in S_j} M_{x,x'} = 1
\end{equation}
for all $i,j \in [c]$.

We will now define a homomorphism matrix for $\K_c \tok \overline{X}$ thus proving the lemma. Define the matrix $H$ entrywise as follows:
\[H_{ix,jx'} = 
\left\{\begin{array}{ll}
M_{x,x'}  & \text{if } x \in S_i, x' \in S_j \\
0 & \text{o.w.}
\end{array}\right.\]
Note that $H$ has $M$ as a principal submatrix, and everything outside of this submatrix is zero. This implies that $H \in \K$ since $H$ can be obtained by tensoring $M$ with $e_1e_1^\intercal \in \cp \subseteq \K$ for $e_1$ the first standard basis vector in $\mathbb{R}^c$ and permuting rows and columns. Considering the matrix $H$ we see that
\[\sum_{x,x'} H_{ix,jx'} = \sum_{x \in S_i, x' \in S_j} M_{x,x}  = 1\]
by equation~(\ref{eqn:M}).

For $i \in [c]$ and $x \ne x'$, we have two cases. In the first case either $x \notin S_i$ or $x' \notin S_i$ and thus $H_{ix,ix'} = 0$. In the second case we have $x,x' \in S_i$ and thus $x \sim x'$, implying that $H_{ix,ix'} = M_{x,x'} = 0$.
 
Suppose that $i \sim j$, i.e.~$i \ne j$, and that $x$ is not adjacent to $x'$ in $\overline{X}$, i.e.~that $x = x'$ or $x \sim x'$ in $X$. To check that $H_{ix,jx'} = 0$, we may assume that $x \in S_i \ne S_j \ni x'$. This implies that $x \ne x'$ and therefore $x \sim x'$. So we have that $H_{ix,jx'} = M_{x,x'} = 0$.

Therefore $H$ is a valid homomorphism matrix for $K_c \tok \overline{X}$ and thus $\alpha^\K(X) = \chi(\overline{X})$ as desired.\qedn

With the above lemma we can quickly prove the main result of this section:

\begin{theorem}\label{thm:stronghom2alpha}
Let $X$ and $Y$ be graphs and $\K$ a frabjous cone. Then $X \tok Y$ if and only if $\alpha^\K(X \ltimes Y) = |V(X)|$.
\end{theorem}
\proof
Recall that $\chi(\overline{X \ltimes Y}) \le |V(X)|$. Therefore, by the above lemma, we have that $\alpha^\K(X \ltimes Y) = |V(X)|$ if and only if $\thetak(X \ltimes Y) = |V(X)|$. Applying Theorem~\ref{thm:hom2theta} completes the proof.\qedn

Though we are not able to prove the analogous theorem to the above for weak conic homomorphisms, we are able to prove one direction.

\begin{theorem}\label{thm:weakhom2alpha}
Let $X$ and $Y$ be graphs and $\K$ a nonnegative frabjous cone. Then $X \tokk Y$ implies $\alpha_\K(X \ltimes Y) = |V(X)|$.
\end{theorem}
\proof
Let $n = |V(X)|$ and let $H$ be a homomorphism matrix for $X \tokk Y$. Since $\alpha_\K(X \ltimes Y) \le n$ always holds, it is necessary and sufficient to construct a homomorphism matrix $H'$ for $K_n \tokk \overline{X\ltimes Y}$. We can assume that the $K_n$ has the same vertex set as $X$. Let
\[H'_{x_1x_2y,x'_1x'_2y'} = 
\left\{\begin{array}{ll}
H_{x_2y,x_2'y'}  & \text{if } x_1 = x_2, x'_1 = x'_2 \\
0 & \text{o.w.}
\end{array}\right.\]
As in the proof of Lemma~\ref{lem:curious}, we have that $H' \in \K$. Also, we have that
\[\sum_{x_2,y,x'_2,y'} H'_{x_1x_2y,x'_1x'_2y'} = \sum_{y,y'} H_{x_1y,x'_1y'} = 1.\]

Now suppose that $x_1 \sim x'_1$ in $K_n$, i.e.~that $x_1 \ne x'_1$, and that $(x_2,y)$ is not adjacent to $(x'_2,y')$ in $\overline{X \ltimes Y}$. To check that $H'_{x_1x_2y,x'_1x'_2y'} = 0$, we may assume that $x_1 = x_2$ and $x'_1 = x'_2$ and thus $x_2 \ne x'_2$. This implies that $(x_2,y) \sim (x'_2,y')$ in $X \ltimes Y$ and therefore $H_{x_2y,x'_2y'} = 0$. This of course implies that $H'_{x_1x_2y,x'_1x'_2y'} = 0$. Therefore $H'$ is a valid homomorphism matrix for $K_n \tokk \overline{X \ltimes Y}$ and we are done.\qedn

Note that the above theorem was proven for $\K = \cpsd$ in~\cite{multiparty}. The above two theorems imply the following corollary which says that in order to check whether $\tok$ and $\tokk$ are equivalent for a given $\K$, it suffices to show that $\alpha^\K = \alpha_\K$.

\begin{corollary}
Let $\K$ be a nonnegative frabjous cone. The relations $\tok$ and $\tokk$ are equivalent if and only if $\alpha^\K(X) = \alpha_\K(X)$ for all graphs $X$.
\end{corollary}
\proof
Obviously, if $\tok$ and $\tokk$ are equivalent, then $\alpha^\K(X) = \alpha_\K(X)$ for any graph $X$ by definition. Conversely, suppose that $\alpha^\K$ and $\alpha_\K$ are the same for all graphs. If $X \tokk Y$, then $\alpha_\K(X \ltimes Y) = |V(X)|$ by Theorem~\ref{thm:weakhom2alpha}. But by assumption we then have that $\alpha^\K(X \ltimes Y) = |V(X)|$, and therefore $X \tok Y$ by Theorem~\ref{thm:stronghom2alpha}. Thus $X \tokk Y$ implies $X \tok Y$ in this case and the converse always holds.\qedn

\section{Discussion}\label{sec:discuss}

We have unified many disparate results concerning homomorphisms, their quantum variations, and other relaxations through the introduction of conic homomorphisms. Furthermore we have explored many properties of conic homomorphisms in general, showing that they behave in many ways similar to the usual homomorphisms. Below we discuss some specific points about conic homomorphisms as well as several important open questions.

\subsection{Restriction on Weak Conic Homomorphisms}\label{sec:weakhomo}

Recall that we only define weak conic homomorphisms for nonnegative frabjous cones. The reason for this is that applying the definition of weak conic homomorphism to a cone which contains matrices with negative entries can yield a relation which is uninteresting/degenerate. To see precisely what this means, consider applying the definition of weak conic homomorphism to the cone $\psd$. Let
\[M = \begin{pmatrix*}[r]
1 & -1 \\
-1 & 1
\end{pmatrix*},
\quad
N = \begin{pmatrix*}[r]
0 & 1 \\
1 & 0
\end{pmatrix*}.\]
Then for large enough $\gamma \in \mathbb{R}$, the matrix
\[H^* = \frac{1}{2}J_n \otimes N + \gamma I_n \otimes M\]
is positive semidefinite and satisfies condition~(\ref{cond:sum}). Moreover, for any graph $X$ with $n$ vertices, the matrix $H^*$ satisfies condition~(\ref{cond:ortho}) for $Y = K_2$. Therefore we have that any graph would have a weak $\psd$-homomorphism to $K_2$. Obviously, this is not desirable. Also, note that $\Theta^\psd = \bartheta$ is not monotone with respect to weak $\psd$-homomorphism, though it is for nonnegative cones as stated in Theorem~\ref{thm:monob}.

To fix this problem, one could add the following constraint to the definition of weak conic homomorphism:
\[H_{xy,xy'} \ge 0 \text{ for all } x \in V(X) \ \& \ y,y' \in V(Y).\]
This would not change the definition for nonnegative cones, since this already holds for those. This extra condition would forbid the use of the matrix $H^*$ above, and allow us to extend Theorem~\ref{thm:monob} to weak conic homomorphism over any frabjous cone. The problem arises with transitivity. To show transitivity for weak conic homomorphism, we would now need to show that the above condition is preserved when we ``compose" conic homomorphisms. If we try to simply modify the proof of Theorem~\ref{thm:trans}, then we would need to show that
\[\sum_{y,y'} H^1_{xy,xy'} H^2_{yz,y'z'} \ge 0\]
for all $x \in V(X)$, and $z,z' \in V(Z)$. It is not clear to the author how to show this.

Note that in the specific case of $\psd$, transitivity does not pose a problem. This is because, if one imposes the additional constraint above, one can show that the weak $\psd$-homomorphism is equivalent to the strong in the same way as for $\dnn$ in the proof of Theorem~\ref{thm:dnnequiv}. However, for general cones this may not be the case.

\subsection{Quantum Correlations}

A (two-party) correlation for input sets $S,T$ and output sets $A,B$ is a joint conditional probability distribution $P(a,b | s, t)$. We say that such a correlation is quantum if there exist positive semidefinite matrices $\rho_{sa}$ for $s \in S$, $a \in A$, and $\sigma_{tb}$ for $t \in T$, $b \in B$, and a unit vector $\psi$ satisfying
\begin{align*}
\sum_{a \in A} \rho_{sa} &= I \text{ for all } s \in S \\
\sum_{b \in B} \sigma_{tb} &= I \text{ for all } t \in T
\end{align*}
and
\[P(a,b|s,t) = \psi^\dagger (\rho_{sa} \otimes \sigma_{tb})\psi.\]
Physically, $P(a',b'|s,t)$ corresponds to the probability that Alice gets outcome $a'$ and Bob gets outcome $b'$ when they respectively perform the quantum measurements $\{\rho_{sa}\}_{a \in A}$ and $\{\sigma_{tb}\}_{b \in B}$ on their shared state $\psi$.

If we consider quantum strategies for the $(X,Y)$-homomorphism game, then we have that $S = T = V(X)$ and $A = B = V(Y)$. It can be shown~\cite{cpsd2qc} that $P(y,y'|x,x')$ is the correlation produced by a perfect quantum strategy for the $(X,Y)$-homomorphism game if and only if the matrix $H$, defined by $H_{xy,x'y'} = P(y,y'|x,x')$, is a homomorphism matrix for $X \tocpsdt Y$. Similarly, $P(y,y'|x,x')$ is a correlation produced by some perfect classical strategy (allowing players access to some shared randomness) if and only if the matrix $H$ defined as above is a homomorphism matrix for $X \tocpt Y$.

If $\K$ is a nonnegative frabjous cone and $H$ is a homomorphism matrix for $X \tok Y$, then $P(y,y'|x,x') := H_{xy,x'y'}$ is a valid joint conditional probability distribution, i.e.~a two-party correlation, although it may not be possible to produce this correlation through the use of classical or even quantum strategies. Interestingly though, this correlation will obey what is known as the non-signalling condition. This condition says that the probability of Alice outputting $y$ upon being given input $x$ is independent of Bob's input. This holds for $P$ defined above because, using the vectors $w$ and $w_{xy}$ for $x \in V(X)$, $y \in V(Y)$ guaranteed by Lemma~\ref{lem:Gram}, we can write the probability of Alice outputting $y$ upon being given input $x$ and Bob being given input $x'$ as
\[\sum_{y'}P(y,y'|x,x') = \sum_{y'} H_{xy,x'y'} = \sum_{y'} \langle w_{xy}, w_{x'y'} \rangle = \langle w_{xy}, w \rangle,\]
which is clearly independent of $x'$.

\subsection{Open Questions}

Perhaps the most important question concerning conic homomorphisms is the following: under what conditions on a nonnegative frabjous cone $\K$ are the relations $\tok$ and $\tokk$ equivalent? A particularly important case of this is whether or not $\tocpsdt$ and $\tocpsdb$ are equivalent. The answer to this question has important implications in zero-error quantum communication. Specifically, an affirmative answer would mean that one can restrict to strategies using projective measurements and maximally entangled states for zero-error communication tasks, and such strategies allow significantly simpler analysis than general quantum strategies.

One approach to this question is to attempt to prove an analog of Theorem~\ref{thm:cp} for completely positive semidefinite matrices. As mentioned above, Laurent and Piovesan~\cite{Piovesan14} showed that the direct analog does not hold in general. However, perhaps it is possible to prove a more restrictive version which would still apply to homomorphism matrices for $X \tocpsdb Y$.\\

Related to this question is the investigation of the completely positive semidefinite cone itself. Since this is such a newly defined object, little is known about it. It seems difficult to even determine if a given matrix belongs to $\cpsd$. It would be useful if one could find alternate descriptions of completely positive semidefinite matrices, similar to the various equivalent characterizations of positive semidefinite matrices. It would also likely be of use to have a description of the extreme rays of $\cpsd$.

There is one property of this cone that is of particular interest: closedness. If $\cpsd$ is closed, then the parameters $\vartheta^\cpsd$ and $\Theta^\cpsd$ are always attained by some solution, which would be of use when trying to determine the value of these parameters for a given graph. Moreover, it has been shown that if $\cpsd$ is closed, then the set of two-party quantum correlations is closed~\cite{cpsd2qc}, and this implies that the well-known Connes' embedding conjecture is true~\cite{connes}.\\

Another interesting question is whether or not $\Thetak$ is multiplicative with respect to the disjunctive and lexicographic products in general. This is true for $\K = \cp$ or $\psd$, but we have not, as of yet, found a way to generalize the proofs for these cases to arbitrary frabjous cones. Even the case of $\K = \dnn$ would be interesting, as this would show that Szegedy's $\barthetap$ is multiplicative.\\

For $\K \in \{\cp, \dnn, \psd\}$, it is known that $\alpha^\K(X) = \lfloor \thetak(X) \rfloor$ for all graphs $X$. This hints that this may be true for all frabjous cones. Further hinting at this is Lemma~\ref{lem:curious}, which seems an odd result if $\alpha^\K$ is not equal to the floor of $\thetak$ in general.\\

Are there other (interesting) frabjous cones? The cones $\dnn$ and $\psd$ give rise to conic homomorphisms which are relaxations of quantum/entanglement-assisted homomorphisms. Can we find other cones which result in relaxations which more closely approximate these relations, but are still tractable? Can we find a frabjous cone which is nested between $\cp$ and $\cpsd$? Many of our results on general conic homomorphisms were inspired by known results for the cases where $\K \in \cones$. Finding more frabjous cones would allow us to apply our results to new variations of graph homomorphism, and may offer new insights to conic homomorphisms in general.

\bibliographystyle{plainurl}
\bibliography{MaxEnt}

\end{document}